\documentclass[11pt]{article}
\RequirePackage[OT1]{fontenc}
\RequirePackage{amsthm,amsmath,natbib}
\RequirePackage[colorlinks]{hyperref}
\usepackage{amsfonts}
\usepackage{amssymb}
\usepackage{amstext}
\usepackage{parskip}
\usepackage{fullpage}

\usepackage{graphicx}
\usepackage{subfig}
\usepackage{multirow}

\usepackage[ruled,linesnumbered,vlined]{algorithm2e}

\usepackage{booktabs}
\usepackage{array}
\usepackage{url}

\providecommand{\keywords}[1]{\textbf{\textit{Keywords---}} #1}

\DeclareFixedFont{\auacc}{OT1}{phv}{m}{n}{12}   
\def\math{${}^\diamond$}
\def\stat{${}^\star$}

\begin{document}
\title{An Extremal Optimization approach to parallel resonance constrained capacitor placement problem
\footnote{Paper published in the 6th IEEE/PES Transmission and Distribution: Latin America, 2012, Montevideo, Uruguay.}}

\author{Andr\'e R. Gon\c{c}alves\math, 
 Celso Cavellucci\stat,\; Christiano Lyra Filho\math, 
Fernando J. Von Zuben\math\;
\\ {\math}School of Electrical and Computer Enginnering, University of Campinas, Brazil
\\{\stat} CLCTEC Consulting, Campinas, Brazil\\
{\footnotesize \{andreric@dca.fee.unicamp.br, celso@clctec.com.br,} \\{\footnotesize chrlyra@densis.fee.unicamp.br, vonzuben@dca.fee.unicamp.br\}} }
\date{}
\maketitle
\begin{abstract}
Installation of capacitors in distribution networks is one 
of the most used procedure to compensate reactive power generated by 
loads and, consequently, to reduce technical losses.
So, the problem consists in identifying the optimal placement and sizing
of capacitors. This problem is known in the literature 
as optimal capacitor placement problem. Neverthless, 
depending on the location and size of the capacitor, it 
may become a harmonic source, allowing capacitor to enter 
into resonance with the distribution network, causing several undesired
side effects. In this work we propose a parsimonious method
to deal with the capacitor placement problem that incorporates resonance
constraints, ensuring that every allocated capacitor will not
act as a harmonic source. This proposed algorithm is based upon a
physical inspired metaheuristic known as Extremal Optimization.
The results achieved showed that this proposal has reached 
significant gains when compared with other proposals that attempt
repair, in a post-optimization stage, already obtained solutions
which violate resonance constraints.
\end{abstract}


\section{Introduction}
\label{sec:intro}

One of the leader causes of technical losses in electric power distribution networks
can be assigned to resistance in distribution and transmission lines, 
subject to active and reactive currents. One possible strategy to reduce 
technical losses is by means of capacitor allocation \citep{Bunch1982}.

Capacitors are used in power distribution networks
with the intention of reducing losses related to reactive power,
power-factor correction, power-flow control and improvement  
of network stability \citep{Madeiro2011}.

The optimal capacitor placement problem (CPP) is to define 
location, size and number of capacitors which should be installed 
in a network, aiming to minimize a function that takes into account
the capacitors installation investment and the reduction in
power loss.

CPP is a non-linear and non-differentiable mixed integer optimization problem with
a set of operating constraints \citep{Yu2004}.
So, traditional optimization methods are not capable of solving large instance
of this hard combinatorial problem. Therefore, techniques conceived to reach
high-quality solutions in a reasonable execution time, such as heuristics
and metaheuristics, emerge as promising approaches.

Some metaheuristics have been successfully applied to CPP, such as:
genetic algorithms \citep{Delfanti2000,Mendes2005}, particle swarm optimization
\citep{Yu2004}, tabu search \citep{Huang1996}, hybrid evolutionary algorithms \citep{Mendes2005}, and genetic fuzzy systems \citep{Das2008}.

However, one of the main drawbacks of population metaheuristics, like genetic algorithms, particle swarm optimization, and ant colony optimization is their high computational cost. Since these 
algorithms handle a population of candidate solutions, a large number of fitness evaluations 
(performance of each candidate solution) must be done, which may be expensive, mainly in real world 
problems \citep{Jin2005}. 

In the CPP, for each fitness evaluation an execution of a power-flow estimation method
is needed. Only after having power-flow estimates we are allowed to figure out the effective
contribution of a given set of allocated capacitors. In real world power distribution 
networks, with thousands of nodes, this procedure becomes computationally expensive.

Another relevant aspect that is commonly relaxed in CPP, or just taken into account in
a post-optimization procedure, is the fact that the capacitor may enter into 
resonance with the power distribution network, depending on some factors, such as the distance
between the capacitor location and the network feeder, and the capacitor size. 
A possible consequence of resonance is an interruption of the power network, that may 
bring several implications to the population and to the electricity concessionaire. So, 
before installing capacitors we need to ensure that the resonance will not show up.

Aiming at mitigating this problem, we propose a parsimonious extremal optimization method 
which takes into account the resonance constraints. Extremal optimization algorithms have
reached good results with relatively small computational cost, when compared with alternative
population metaheuristics \citep{Lu2009} and \citep{Lu2007}.

The remainder of the text is structured as follows: in the next section we discuss about
resonance in distribution network. The mathematical modeling of the resonance constrained
CPP is described in Section~\ref{sec:modelagem}. Extremal Optimization is depicted in 
Section~\ref{sec:oe}. Our proposal is outlined in Section~\ref{sec:algoritmo_proposto}. 
Experiments carried out are described in Section~\ref{sec:experimentos}, and the results are analyzed in 
Section~\ref{sec:resultados}. Finally, concluding remarks and future works are pointed out
in Section~\ref{sec:conclusao}.

\section{Resonance in power distribution networks}
\label{sec:ressonancia}

The size and location of capacitors are critical factors in a distribution system's
response to the harmonic sources \citep{Gonen1986}. Combination of capacitors and
the system reactance cause parallel resonance frequencies 
for the circuit. The possibility of resonance between a capacitors and the 
rest of the system, at a harmonic frequency, can be determined by calculating 
equal order of harmonic $h$ at which resonance may take place, given by \citep{Gonen1986}

\begin{equation}
	h = \sqrt{\frac{S_{cc}}{Q_c}},
\label{eq:harmonica}
\end{equation}

\noindent where $S_{cc}$ is the three-phase short-circuit power of system at 
the place of capacitor installation, in VA, and $Q_c$ is the capacitor size, given in var.

Parallel resonance frequency, $f_p$, can be expressed as

\begin{equation}
	f_p = f_1\cdot h,	
\label{eq:freq_harm}
\end{equation}

\noindent where $f_1$ is the fundamental frequency (60 Hz in Brazil). Replacing  Eq.~\ref{eq:harmonica}
in Eq.~\ref{eq:freq_harm}, we have

\begin{equation}
f_p = f_1\cdot \sqrt{\frac{S_{cc}}{Q_c}}.
\label{eq:freq_harm_2}
\end{equation}

The characteristics of the most common harmonic loads in power distribution networks 
include the third, fifth and seventh harmonics. Thus, the purpose is to avoid 
these \emph{odd} harmonics. In this work, we consider a range of 10 Hz above and
below the frequency $f_p$, that is, the $n$-th harmonic is characterized by the interval, 
[$(n\cdot f_1)-10 Hz, (n\cdot f_1)+10 Hz$].


Problems resulting from harmonics include (among others) \citep{Mahmoud1983}:

\begin{itemize}
	\item Inductive interference with telecommunication systems;
	\item Capacitor failure from dieletric breakdown or reactive power overload;
	\item Excessive losses in - and heating of - induction and synchronous machines;
	\item Dieletric instability of insulated cables resulting from harmonic overvoltages
		on the system.
\end{itemize}

Effects of harmonics in capacitors include: ($i$) capacitor overheating, ($ii$)
overvoltage in the capacitor, and ($iii$) losses in the capacitor.

Among the harmonics control techniques we can cite \citep{Gonen1986}: (1) strategically identify 
the location of the capacitors for intallation, (2) select the capacitors size properly, and 
(3) remove capacitors that act as harmonic sources. Our proposal include
the first two control techniques, avoiding the third one. As we shall see in 
Section~\ref{sec:resultados}, the third strategy is not always a good choice.

\section{Mathematical modeling}
\label{sec:modelagem}

Capacitor placement problem seeks to minimize the cost of power losses 
and the investment made in capacitors. We can note that both objectives
are specified in monetary terms. In fact, the higher the investment in
capacitors, the smaller the power loss and, consequently, the smaller the 
cost of power loss. So, the algorithm needs to find a balanced solution between
these two criteria: cost of losses and investment made in capacitors.
It is worth mentioning that although this problem could be seen from a multiobjective 
perspective, we performed an equal-weighted criteria approach.


The entire network loss can be taken from Eq.~\ref{eq:perdas}.

\begin{equation}
f_p(P,Q,V) = \sum_{n \in N} \sum_{a \in A_n} r_a\left( \frac{P_{(n,a)}^2 + Q_{(n,a)}^2}{V_{(n,a)}^2}\right ),
\label{eq:perdas}
\end{equation}

\noindent where $N$ is the number of graph nodes (since the distribution network is represented
by a graph \citep{Cavellucci1997}), $A_n$ correspond to the arc set which emanates from node $n$, $r_a$ is the resistance 
in the path $a$, $P_{(n,a)}$ and $Q_{(n,a)}$ are active and reactive power, respectively,
flowing through arc $a$, for a given period of time (1 hour). For the sake of simplicity,
we will assume that voltage values ($V_{(n,a)}$) are approximately 1 pu, for all network nodes.
Now, reformulating Eq.~\ref{eq:perdas} to include capacitors, the network losses are defined as


\begin{equation}
 f_p(P,Q,\bar{Q}) = \sum_{n \in N} \sum_{a \in A_n} r_a (P_{(n,a)}^2 + Q_{(n,a)}^2 - \bar{Q}_{(n,a)}^2).
\label{eq:fp}
\end{equation}

In this work, the cost of the entire network loss is calculated for a period of
one year (8760 hours). From this, we are able to define the optimization
problem to be solved, which can be described by Eq.~\ref{eq:modelo_matematico}.

\begin{equation}
\begin{aligned}
& \underset{u_n}{\text{minimize}}
& & \left \{8.76\cdot Cost\cdot f_p(P,Q,\bar{Q}) + \sum_{n\in N} f_c(u_n) \right  \} \\
& \text{subject to}
&& P_{n-1}=P_{L_n} + \sum_{a \in A_n} P_{a} \\
&&& Q_{n-1}=Q_{L_n} - \bar{Q}_{n} + \sum_{a \in A_n} Q_{a} \\
&&& V_a \approx \textrm{1 p.u.} \\
&&& \bar{Q} \in \Omega_{\bar{Q}} \\
&&& h_i \text{ is even}, \forall i=1,2,...,N \\
\end{aligned}
\label{eq:modelo_matematico}
\end{equation}

\noindent where $\bar{Q}_{n}$ is the installed capacitor at the location $n$, {\it Cost} is the energy price
(in U\$) per MWh, $h_i$ is the harmonic on the $i$-th network bar, that will be better discussed in Section~\ref{sec:algoritmo_proposto},
and $f_c(u)$ is amortised capacitor cost. In the problem formulation, Eq.~\ref{eq:modelo_matematico}, the 8.76 value concerns
the number of hours in a year (8760) divided by 1000, since losses are measured
in kWh and the energy price in U\$/MWh.

The amortised capacitor cost is given by:
 
\begin{equation}
f_c(u_n) = \begin{cases}
\frac{i\cdot g_c(u_n)}{1-1/(1+i)^k},&\text{if there is a capacitor in }n\\ 
 0, & \text{otherwise}
\label{eq:custo_anual}
\end{cases}
\end{equation}

The capacitors cost described in Eq.~\ref{eq:modelo_matematico} corresponds to the 
total cost subtracted by annual gain. So, it is necessary to define an amortization 
constant $k$ for the equipment and an interest annual rate $i$ (Eq.~\ref{eq:custo_anual}). The period of recovery
generally corresponds to the useful life of the equipment. In our work we have adopted used an annual rate 
equal to 0.12 (12\%) and a period $k$ equal to five years.

For the cost of the capacitors $g_c(u_n)$ it was used the following table
of capacitor available for installation ($\Omega_{\bar{Q}}$ set): 

\begin{table}[!htb]
\caption{Types of capacitor available for installation.}
\label{tab:capacitores}
\centering
\begin{tabular}{cccc}\toprule
   Type (u$_n$) & Size (kvar)     & Cost (U\$) & Cost/kvar\\ \toprule  
      1 	& 150             & 1498     & 10,00 \\
      2 	& 300             & 1604     & 5,35 \\
      3 	& 450             & 1620     & 3,60 \\
      4 	& 600             & 1823     & 3,04 \\
      5 	& 900             & 2550     & 2,83 \\
      6   & 1200            & 2955     & 2,46 \\            
  \bottomrule
\end{tabular}
\end{table}

The resonance constraint is incorporated into the problem by means of a
frequency scanning technique, looking for odd frequencies, given that they are commonly
found in distribution networks. This procedure is better explained 
in Section~\ref{sec:algoritmo_proposto}.

\section{Extremal Optimization}
\label{sec:oe}

Extremal optimization (EO) \citep{Boettcher1999,Boettcher2000} is a general-purpose local search 
heuristic based upon recent progress in understanding far-from-equilibrium phenomena in terms of 
self-organized criticality (SOC) \citep{Boettcher2002}. The dynamic of EO was inspired 
by self-organized criticality, a concept introduced to describe emergent complexity in physical
systems, where an optimized structure emerges naturally by simple selection against the extremelly bad.
EO method, as well as {\it Simulated Annealing} (SA) \citep{Kirkpatrick1983} and {\it Genetic Algorithm}
(GA) \citep{Goldberg1989}, are inspired by observations of natural systems.

Unlike GA, which is a population algorithm, EO handles only one solution at a time and 
seeks to improve the quality of this solution through local perturbations. Originally,
this algorithm was proposed to deal with combinatorial optimization problems, particularly
problems which can be represented by a graph. In these applications, EO has been shown
competitive with more elaborate general-purpose heuristics on testbeds of constrained
optimization problems with up to $10^5$ variables, such as bipartitioning, coloring,
and satisfiability \citep{Boettcher2002}.

In a graph representation of CPP, variables are nodes and the influence between variables 
are represented by the arcs. So, a node perturbation will affect directly its neighbors
(parent and children in the case of a tree).

In evolutionary algorithms, a quality measure is assigned to each solution, called fitness. 
Differently from these approaches, EO assigns a fitness to each variable (although that is 
not essential \citep{Boettcher2000}), $\lambda_i$, being the \emph{total cost} of solution, 
$C(S)$, obtained as follows

\begin{equation}
C(S) = \sum_{i=1}^{n} \lambda_i.
\end{equation}

As aforementioned, perturbations are made in a selected variable, the one with the smallest fitness 
value (maximization problem). Due to its influence on the neighboring variables (neighbor graph nodes),
this perturbation will also reflect in these variables. The pseudocode of the basic EO algorithm is 
presented in Algorithm~\ref{alg:eo}.\\

	\begin{algorithm}
		\DontPrintSemicolon
		\KwResult{$S_{best}$ and $C(S_{best})$.} 
		\Begin{
			Define an initial solution $S$; and set $S_{best}$ = $S$\\
			\Repeat{stopping condition met}{
				Evaluate $\lambda_i$ for each variable $x_i$ \\
				Find $j$ satisfying $\lambda_j \leq \lambda_i$, for all $i$\\
				Choose $S'\in N_S$ such that $x_j$ must change\\
				$S \leftarrow S'$ \\
				\If{ $C(S) < C(S_{best})$ }{
					$S_{best} \leftarrow S$\\
				}
			}
		}
		\caption{Framework of Extremal Optimization.}
		\label{alg:eo}
	\end{algorithm}

An EO variation was proposed by \citep{Boettcher2001}, such that the selection of the variable
with the worst fitness is not done determiniscally, but follows a probability distribution,
favoring the worst variables, but also providing a chance to the others. A 
power law probability distribution was employed to perform this task, given by

\begin{equation}
P(k)\propto k^{-\tau}\;\;\;(1\leq k \leq n),
\label{eq:leipot}
\end{equation}

\noindent where $\tau \in \Re$ is a probability distribution parameter. This 
variation was called $\tau$-EO. For $\tau=0$, the algorithm is a local-random search 
algorithm. On the other hand, for $\tau\rightarrow\infty$, this is a deterministic local search
algorithm, where the variable with the worst fitness is always selected to be updated.
For the traveling salesman problem, $\tau$ values between 1.6 and 2 achieved good
results \citep{Boettcher2002}.

Another variation of EO, proposed by \citep{Chen2007}, was adopted here.
In this method, besides employing a probability distribution to select a variable
to be updated, another power law distribution is applied to select, among a set of 
neighboring candidate solutions, which one will replace the current solution.
This approach lead to a more informative EO version and, as a result, it is possible to 
find good solutions faster than the original version.

\section{EO approach to resonance constrained CPP}
\label{sec:algoritmo_proposto}

In this section, we describe an extremal optimization algorithm to deal 
with resonance constrained CPP. The pseudocode of this proposal is shown in 
Algorithm~\ref{alg:proposta}.

The main reason to employ an extremal optimization algorithm, instead of any 
other metaheuristics, is due to the fact that EO was initially developed
to large combinatorial optimization problems which can be represent by a graph,
that is precisely the inherent nature of CPP. Another reason, which will become clear later, 
is the reduction in computational resource when handling resonance constraints. 

\begin{algorithm}
		\KwResult{$S_{best}$ and $C(S_{best})$} 
		\Begin{
		    Set the initial solution, $S$, as a null vector\\
		   $S_{best} \leftarrow S$\\
		   Calculate fitness, $\lambda_i$, of each node (variable)\\
		    \Repeat{stopping criterion is not satisfied}{
		    	Sort variables by increasing fitness\\
		    	Select a variable by following a power law pdf\\
		    	Generate neighbors from the selected variable\\
		    	Calculate $\lambda_i$ and \emph{total cost} of all neighbors\\
		    	Sort neighbors by decreasing fitness\\
		    	Select a neighbor by following a power law pdf\\
		    	Update $S_{best}$ if needed\\
		    }
		}
        \caption{Proposed algorithm.}
        \label{alg:proposta}
\end{algorithm}

Solutions are represented by an integer-valued vectors belonging to the interval [0,6], 
being 0 if there is no capacitor in that network location and 1 to 6 if it has one
of those six types of capacitor shown in Table~\ref{tab:capacitores}.

For the initial solution we assume that there is no capacitor installed, that is, the initial 
solution is an $n$-dimensional vector consisting of all zeros, where $n$ is the number of nodes in 
the tree, representing the distribution network.

\begin{figure}[!ht]
	\centering
	\includegraphics[scale=0.35]{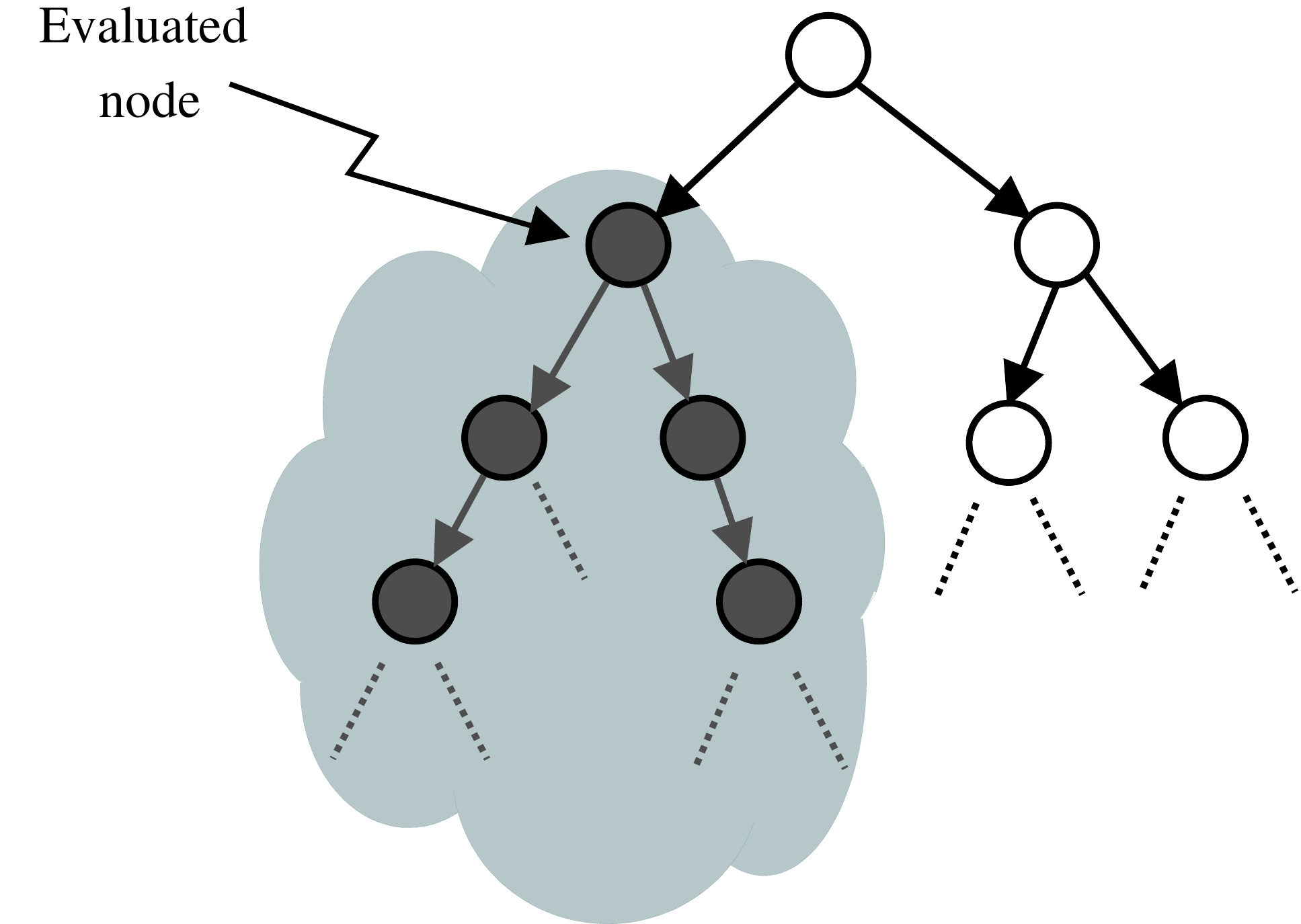}
	\caption{Subnetwork used to compute the node loss.}
	\label{fig:fitness_variable}
\end{figure}

For each node of the tree, the cost function is evaluated. The losses $f_p(\cdot)$ of a given 
node are accumulated up to this network point, which is calculated by a power-flow estimation algorithm.
Figure~\ref{fig:fitness_variable} shows how the losses of a node is calculated. The idea is to isolate the 
subnetwork rooted by the interested node and calculate the losses for this subnetwork.

In this work, a simplified version of the power-flow estimation algorithm proposed by Baran and Wu 
\citep{Baran1989} was used.

After the cost function be evaluated for each node, they are sorted by increasing fitness and 
one node is chosen, following a power law probability distribution (Eq.~\ref{eq:leipot}), to 
be perturbated. Perturbation mechanism is viewed as a local exploration of the current solution, 
generating a set of neighboring candidate solutions $N_S$, which are slightly different from $S$.

At this point the resonance control procedure is applied. From the selected node, 
its neighborhood is generated by Algorithm~\ref{alg:generate_neighbors}.

\begin{algorithm}
		\DontPrintSemicolon
		\KwIn{$i, S$} 
		\KwResult{$N_S$.} 
		\Begin{		
		    $N_S \leftarrow \emptyset$ \\
		    \If(\tcc*[f]{remove capacitor}){ $S[i] \neq 0$ }{
		    	$S'[i] \leftarrow 0$\\
		    	$N_S \leftarrow N_S \cup S'$\\
		    }
		    \If(\tcc*[f]{increase capacitor size}){$S[i] \neq 0$}{
		    	$S'[i] \leftarrow min(S[i]+1,6)$\\
		    	$N_S \leftarrow N_S \cup S'$\\
		    }
		    \If(\tcc*[f]{decrease capacitor size}){$S[i] \neq 0$}{
		    	$S'[i] \leftarrow max(S[i]-1,0)$\\
		    	$N_S \leftarrow N_S \cup S'$\\
		    }
		    \If(\tcc*[f]{install a new capacitor}){ $S[i] == 0$}{
		    	$S'[i] \leftarrow \texttt{rand\_integer(1,6)}$\\
		    	$N_S \leftarrow N_S \cup S'$\\
		    }
		    \If(\tcc*[f]{shift the capacitor to its parent}){ $S[i] \neq 0$}{
		    	$S'[parent_i] \leftarrow S[i]$\\
		    	$S'[i] \leftarrow 0$\\
		    	$N_S \leftarrow N_S \cup S'$\\
		    }
		    \If(\tcc*[f]{shift the capacitor to its children}){$S[i] \neq 0$}{
		   	\For{$k$ {\bf in} children}{
			    	$S'[k] \leftarrow S[i]$\\
		    		$S'[i] \leftarrow 0$\\
			    	$N_S \leftarrow N_S \cup S'$\\
			}}
    }
    \caption{Generate neighbors.}
    \label{alg:generate_neighbors}
\end{algorithm}

It can be seen that there is a minimum and maximum quantity of solutions in 
its neighborhood. When there is no capacitor allocated in that position, just
one neighbor is generated (lines 14-17). Otherwise, we get the maximum number of 
neighbors, 4+$k$, where $k$ is the number of descendant nodes (children). 

A solution is only inserted in the neighbor set if, and only if, the parallel resonance constraint
is satisfied. That is a penalty function constraint-handling method, which is known as 
\emph{death-penalty} approach \citep{MezuraMontes2011}. Even though there are many
methods to handle constraints in nature-inspired algorithms, death-penalty
is a very simple one and, as we will see in the results section, it has reached good performance.
Algorithm~\ref{alg:res_constraint} shows the routine developed to check the resonance constraint.

\begin{algorithm}
		\KwIn{$n$} 	   
		\KwResult{\emph{satisfy}} 
		\Begin{		
		    $f_p = f_1\cdot \sqrt{S_{cc}(n)/\bar{Q}(n)}$\\
		    $h$ = \texttt{round(}$f_p$/$f_1$\texttt{)}\\
		    \If{$h$ is even}{
		    	\emph{satisfy} = True\\
		    }{
		    	\emph{satisfy} = False\\
		    }
		}
    \caption{Check resonance constraint routine.}
    \label{alg:res_constraint}
\end{algorithm}


A key point of our extremal optimization approach is that only one network node has to be 
verified at each generation, unlike a genetic algorithm after the application of a uniform crossover operator, 
for example, in which a high number of nodes must be checked. Thus, the analysis over resonance constraint 
is made in only one network point whatever the network size. This may save 
a substantial computational resource, mainly in real power distribution networks,
composed of thousands of nodes. That is a significant advantage of extremal optimization 
approches in relation to another metaheuristic.

Once the set of neighboring solutions has been determined, they are evaluated and sorted 
by decreasing \emph{total cost $C(\cdot)$} and, then, one of these neighbors 
is selected to replace the current solution $S$, following another power law 
distribution, defined as \citep{Zeng2010}

\begin{equation}
	P(k) \propto e^{-\mu k}\;\;\;(1\leq k \leq n),
	\label{eq:leipotexp}
\end{equation}

\noindent where $\mu \in \Re_+$ is a distribution parameter. 

Update the best solution found so far, $S_{best}$, if necessary, and check the 
stopping criterion.

\section{Experiments}
\label{sec:experimentos}

\begin{figure}[!ht]
	\centering
	\includegraphics[scale=0.4]{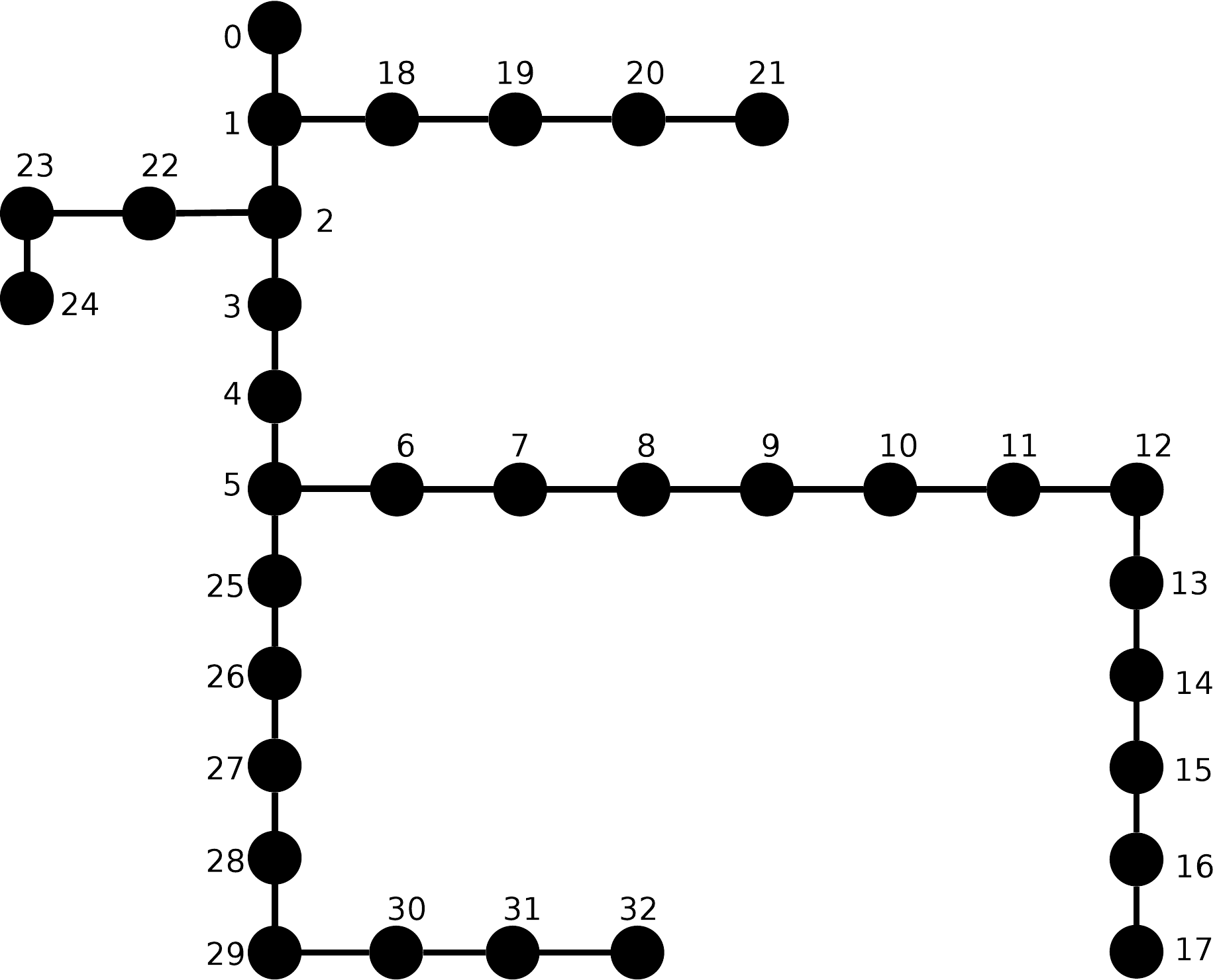}
	\caption{Power distribution network described by Baran and Wu \citep{Baran1989}.}
	\label{fig:redeWu}
\end{figure}

Our algorithm was applied to a power distribution network initially described in Baran and Wu \citep{Baran1989}.
This network is composed of 33 nodes and 34 arcs. Although of small size, that network was widely used in
the literature to compare performance of algorithms designed to cope with power distribution system problems, like CPP 
and network reconfiguration \citep{Madeiro2011,Jeon2002,Jeon2004}. Baran and Wu (BW) network is illustrated in Figure~\ref{fig:redeWu}.

The performance of the considered algorithms was analyzed in relation to the number of function evaluations (FEs),
instead of number of generations. This is due to the fact that populational algorithms, like GA, execute
a large amount of function evaluations per generation, depending on the population size, whereas single individual
algorithms, such as EO, perform a reduced amount of function evaluations. Function evaluation phase is
very expensive in this type of problem, since it is necessary to run a power-flow estimation procedure for each
candidate solution.

\subsection{Comparative analysis}
\label{sec:algoritmoscomparacao}

To comparativelly analyze the performance of our EO method, a memetic 
algorithm proposed by \citep{Mendes2005}, designed to deal with CPP, was implemented.
Despite not taking the resonance constraint into account, the memetic algorithm has reached
good results for real power distribution networks. To get a final feasible solution,
some post-optimization strategies (described in what follows) are applied, if the final
solution does not satisfies the resonance constraints.

This memetic algorithm is a genetic algorithm that uses a hierarchically
tree-structured population, composed of 13 individuals and a local search
procedure, applied to the best individual in the population (placed at the 
tree root). Uniform crossover and punctual mutation operators are employed.

Some of the commonly used strategies (in practice) to repair a final unfeasible solution
for CPP is described below. These strategies was named STRTG1, STRTG2 and STRTG3, respectively.

\begin{itemize}
\item {\bf STRTG1:} Removing capacitors that has entered into resonance with the distribution network;
\item {\bf STRTG2:} Shifting capacitors that has entered into resonance with the network to their respective
\emph{parents};
\item {\bf STRTG3:} Shifting capacitors that has entered into resonance with the network to their respective
\emph{children}.
\end{itemize}

If the constraint is still not satisfied, these solutions are dropped from analysis.
Otherwise, their gains will be compared with the one achieved by our algorithm.

It is worth mentioning that STRTG1 always will return a feasible solution, whereas in 
the remaining two it is not guaranteed.

Both algorithms and the power-flow estimation procedure were implemented in Python 2.7 using Numpy. 
Simulations were performed on a Intel Core tm 2 Quad Q6600 @ 2,40 Ghz and 2 GB RAM.

\section{Results and Discussion}
\label{sec:resultados}

In our simulations the following parameter values were used for GA (the 
same values used by Mendes {\it et al.} \citep{Mendes2005}): $rate_{cross}=1.5$, 
$p_{mut}=0.1$ and structured population with 13 individuals. In the case 
of EO, the values (defined by a grid search procedure) are $\tau=2$ and $\mu=0.5$.
For both algorithms, the number of fitness evaluations was limited in 50.000
and the results were analyzed over 30 independent runs.

Table~\ref{tab:results} shows the results in terms of mean and standard deviation of the
amount of saved money when using the best solution reached by each algorithm, 
varying linearly the energy price, from 50 to 150 dollars per MWh. The absence of results
for MA+STRTG3 is because this approach was not able to produce feasible solutions.

\begin{table*}[!htb]
\caption{Mean and standard deviation of the amount of saved money for each algorithm 
when varying the energy price.}
\label{tab:results}
\centering
\footnotesize
\begin{tabular}{cccccc}\toprule
  {\bf  Energy}  & \multirow{2}{*}{\bf EO} & \multirow{2}{*}{\bf MA+STRTG1} & \multirow{2}{*}{\bf MA+STRTG2} & \multirow{2}{*}{\bf MA+STRTG3} \\
   {\bf Price (U\$)} &  &  & \\ 
   \toprule  
   50  & {\bf 22,573.80} ($\pm$ {\bf 229.01}) & 17,643.85 ($\pm$ 9.87)    & 19,902.41 ($\pm$ 9.87)  & --  \\ 
  60  & {\bf 27,591.14} ($\pm$ {\bf 268.27}) & 21,424.15 ($\pm$ 7.2e-12) & 24,217.54 ($\pm$ 1.4e-11) & -- \\ 
  70  & {\bf 32,547.35} ($\pm$ {\bf 305.18}) & 25,200.06 ($\pm$ 13.82)   & 28,528.27 ($\pm$ 13.82)   & -- \\ 
  80  & {\bf 37,419.25} ($\pm$ {\bf 289.14}) & 28,967.14 ($\pm$ 19.74)   & 32,830.17 ($\pm$ 19.74)   & -- \\ 
  90  & {\bf 42,354.99} ($\pm$ {\bf 345.82}) & 32,770.64 ($\pm$ 4.74) & 37,168.50 ($\pm$ 4.74) & -- \\ 
 100  & {\bf 47,536.34} ($\pm$ {\bf 364.21}) & 36,597.07 ($\pm$ 6.37) & 41,529.75 ($\pm$ 6.37) & -- \\ 
 110  & {\bf 52,332.52} ($\pm$ {\bf 425.12}) & 40,421.99 ($\pm$ 5.90) & 45,889.49 ($\pm$ 5.90) & -- \\ 
 120  & {\bf 57,468.45} ($\pm$ {\bf 356.68}) & 44,248.70 ($\pm$ 6.19) & 50,251.03 ($\pm$ 6.19) & -- \\ 
 130  & {\bf 62,285.08} ($\pm$ {\bf 413.51}) & 48,074.14 ($\pm$ 6.85) & 54,611.30 ($\pm$ 6.85) & -- \\ 
 140  & {\bf 67,365.68} ($\pm$ {\bf 473.19}) & 51,899.50 ($\pm$ 7.50) & 58,971.48 ($\pm$ 7.50) & -- \\ 
 150  & {\bf 72,301.80} ($\pm$ {\bf 554.78}) & 55,728.08 ($\pm$ 6.94) & 63,334.88 ($\pm$ 6.94) & -- \\ 
  \bottomrule
\end{tabular}
\end{table*}


It can be seen that our proposal has reached better feasible solutions than
memetic algorithm with post-optimization repair procedure, for all problem
configurations. The difference in performance between algorithms becomes higher when
energy price increases. 

A statistical comparison between these results was done by means of the {\it t-test}.
For all energy prices, the proposed EO-based algorithm statistically outperformed
the other ones with {\it p-value} around 1e-50, i.e., the equality hypothesis is
utterly rejected.

The obtained results gives us some evidences that STRTG3 is not a good post-optimization
repair strategy, whereas STRTG1 and STRTG2 have succeeded in getting feasible solutions
in all cases. Among the strategies 1 and 2, shifting capacitors to its parents (STRTG2)
seems to be the more promising one, once it has reached better solutions for all situations.

\begin{figure*}[htb]
  \centering
  \subfloat[Proposed EO approach.]{\label{fig:eo_qte_varcost}\includegraphics[width=0.45\textwidth]{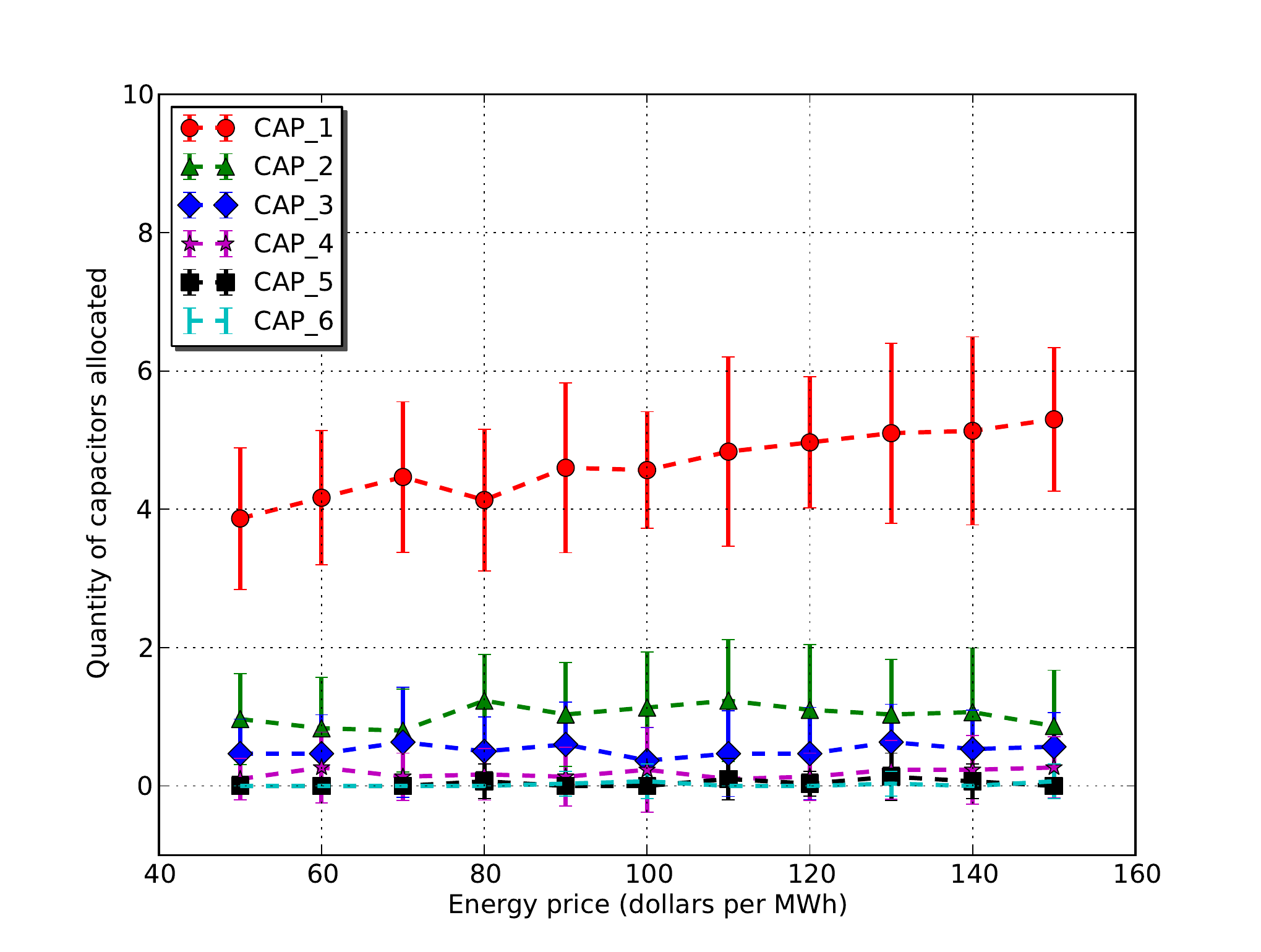}}
  ~ 
  \subfloat[MA+STRTG1.]{\label{fig:stg1_qte_varcost}\includegraphics[width=0.45\textwidth]{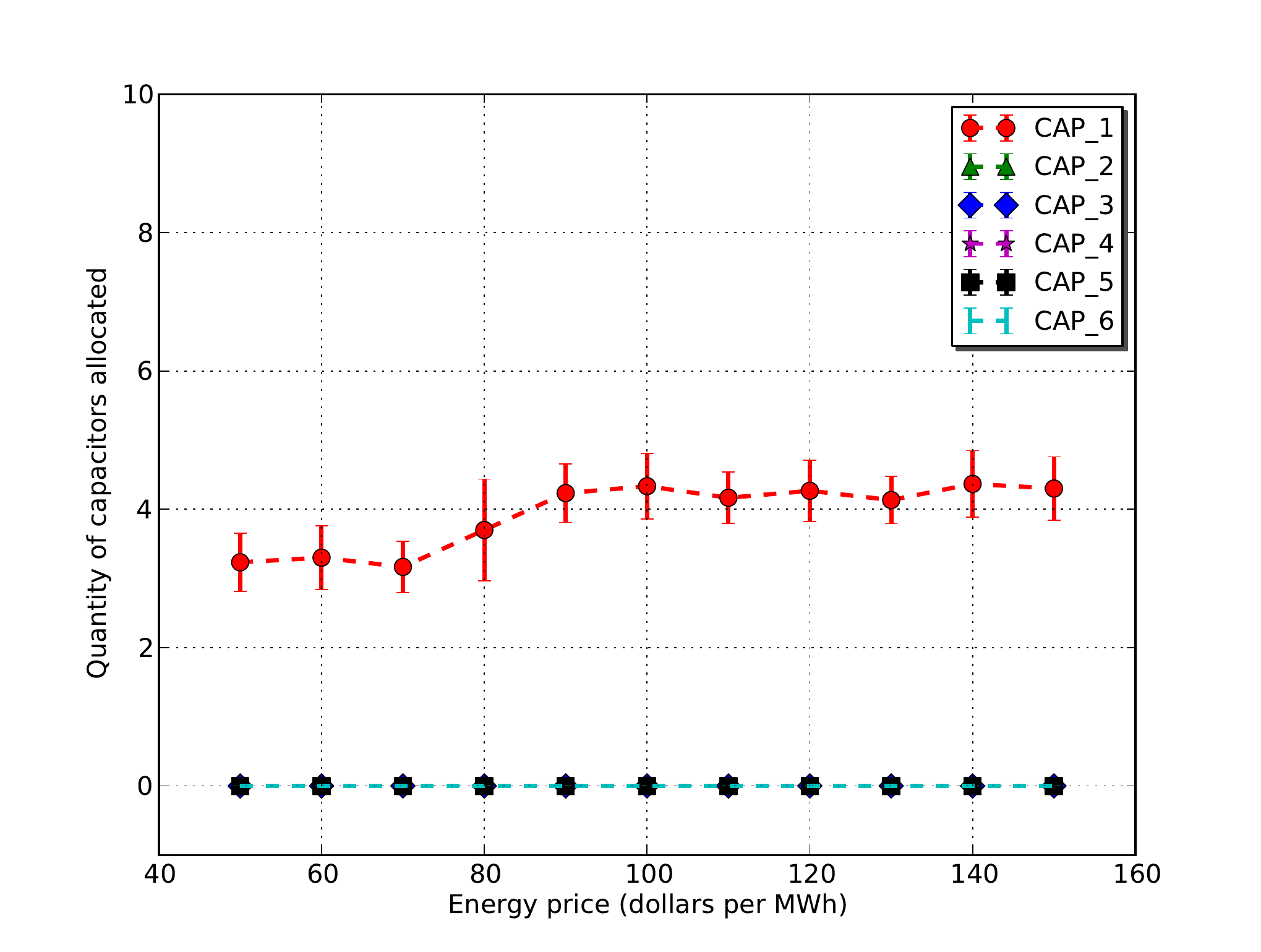}}

  \subfloat[MA+STRTG2.]{\label{fig:stg2_qte_varcost}\includegraphics[width=0.45\textwidth]{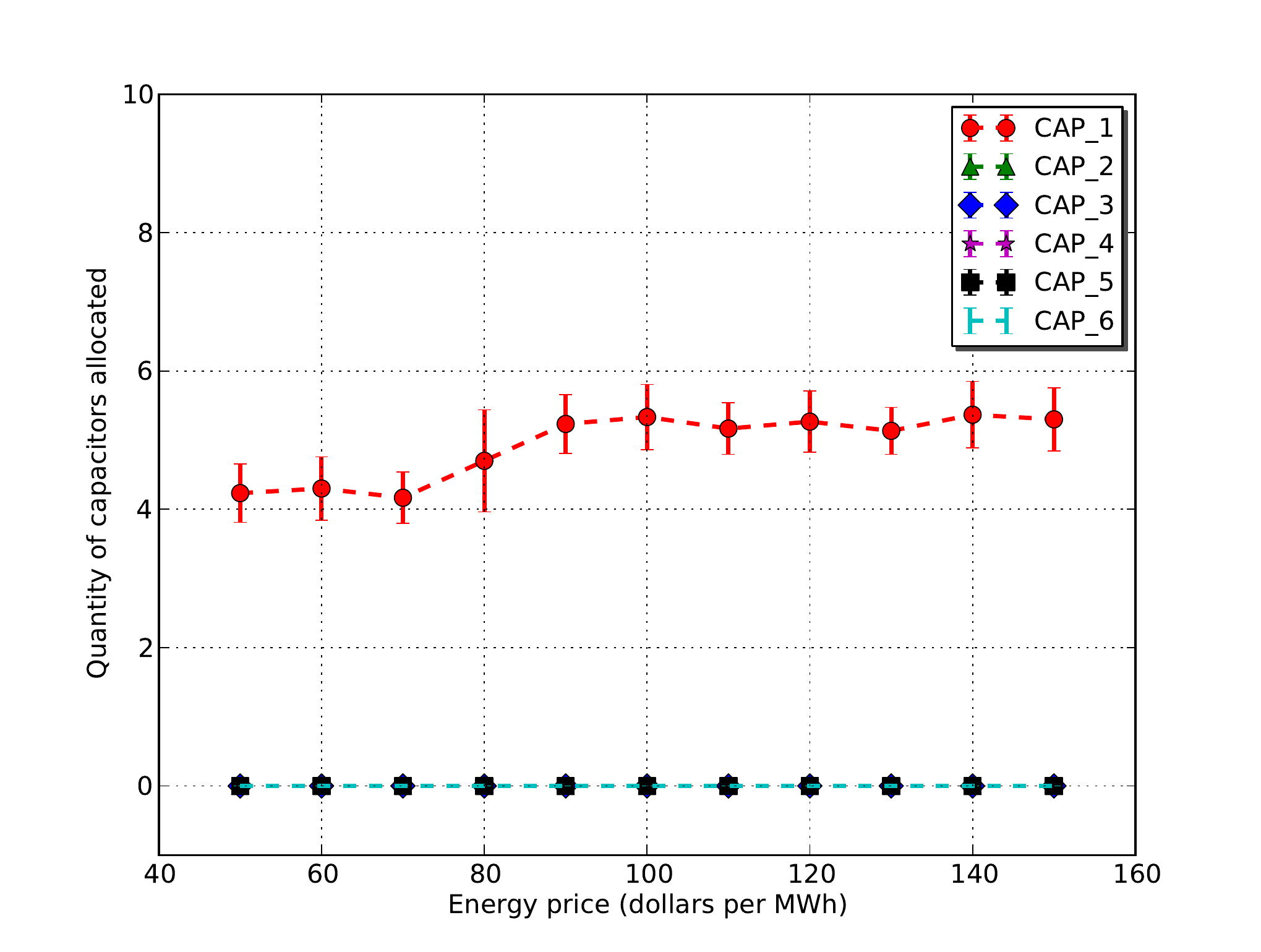}}
  \caption{Quantity of each type of capacitors allocated by the algorithms varying the energy price.}
  \label{fig:type_caps_allocated}
\end{figure*}

Figures (\ref{fig:eo_qte_varcost}), (\ref{fig:stg1_qte_varcost}), and (\ref{fig:stg2_qte_varcost})
depicts the mean and standard deviation of the quantity of each capacitor types allocated by
each algorithm, varying the energy price. Due to the fact that both, STRTG1 and STRTG2, are 
slight modification of the same (unfeasible) solution, their curves have a quite similar behavior.

As we can see in these plots, the number of capacitors installed by EO is higher, in general, than the number proposed by 
the three other approaches. In addition, its solutions are composed of a larger variety of capacitor types in relation 
to the solutions of the other approaches. It can might be the result of a better exploration of the search space.

\begin{figure*}[htb]
  \centering
  \subfloat[Proposed EO approach.]{\label{fig:eo_qte_varcaps}\includegraphics[width=0.45\textwidth]{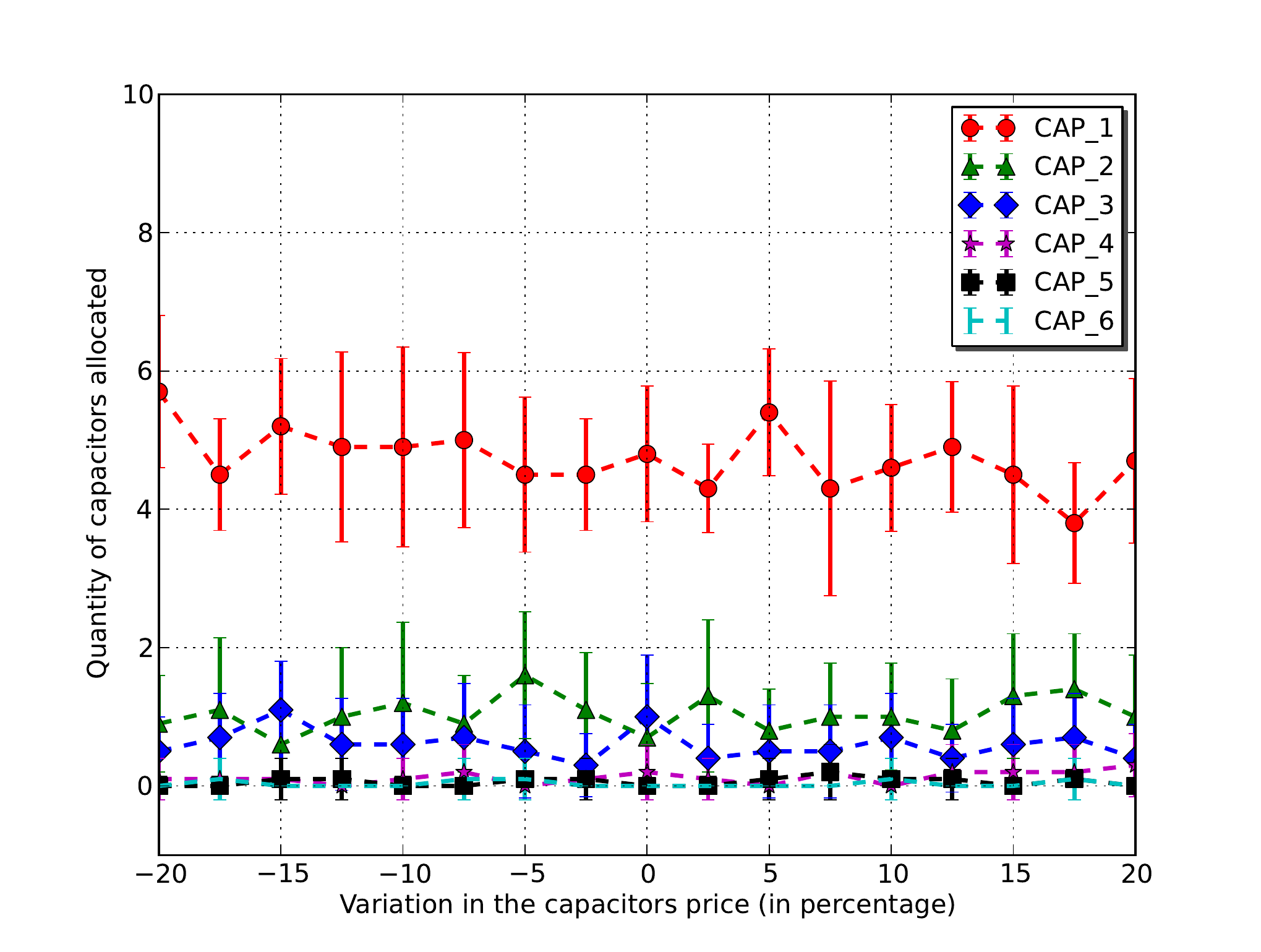}}
  ~ 
  \subfloat[MA+STRTG1.]{\label{fig:stg1_qte_varcaps}\includegraphics[width=0.45\textwidth]{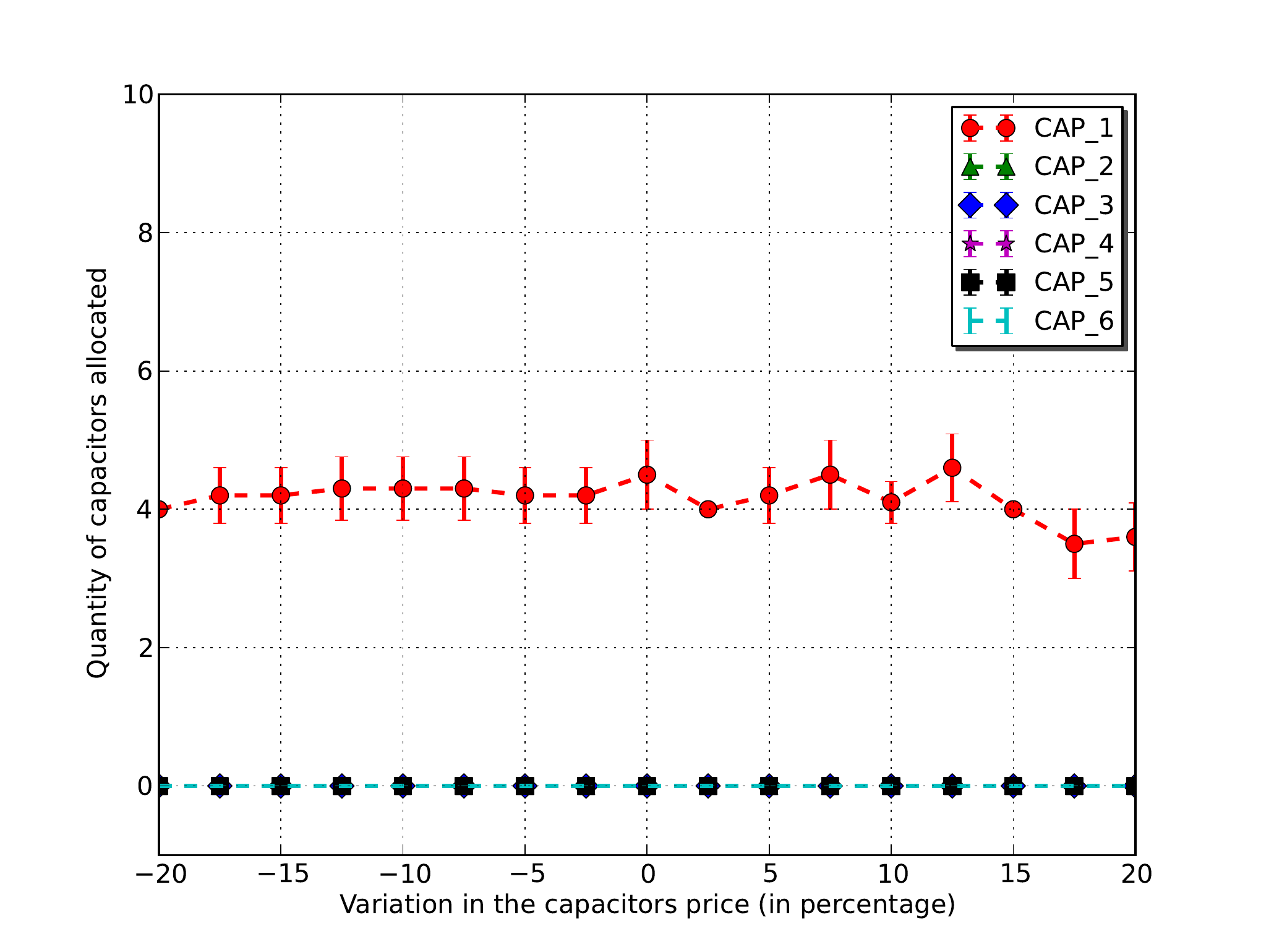}}

  \subfloat[MA+STRTG2.]{\label{fig:stg2_qte_varcaps}\includegraphics[width=0.45\textwidth]{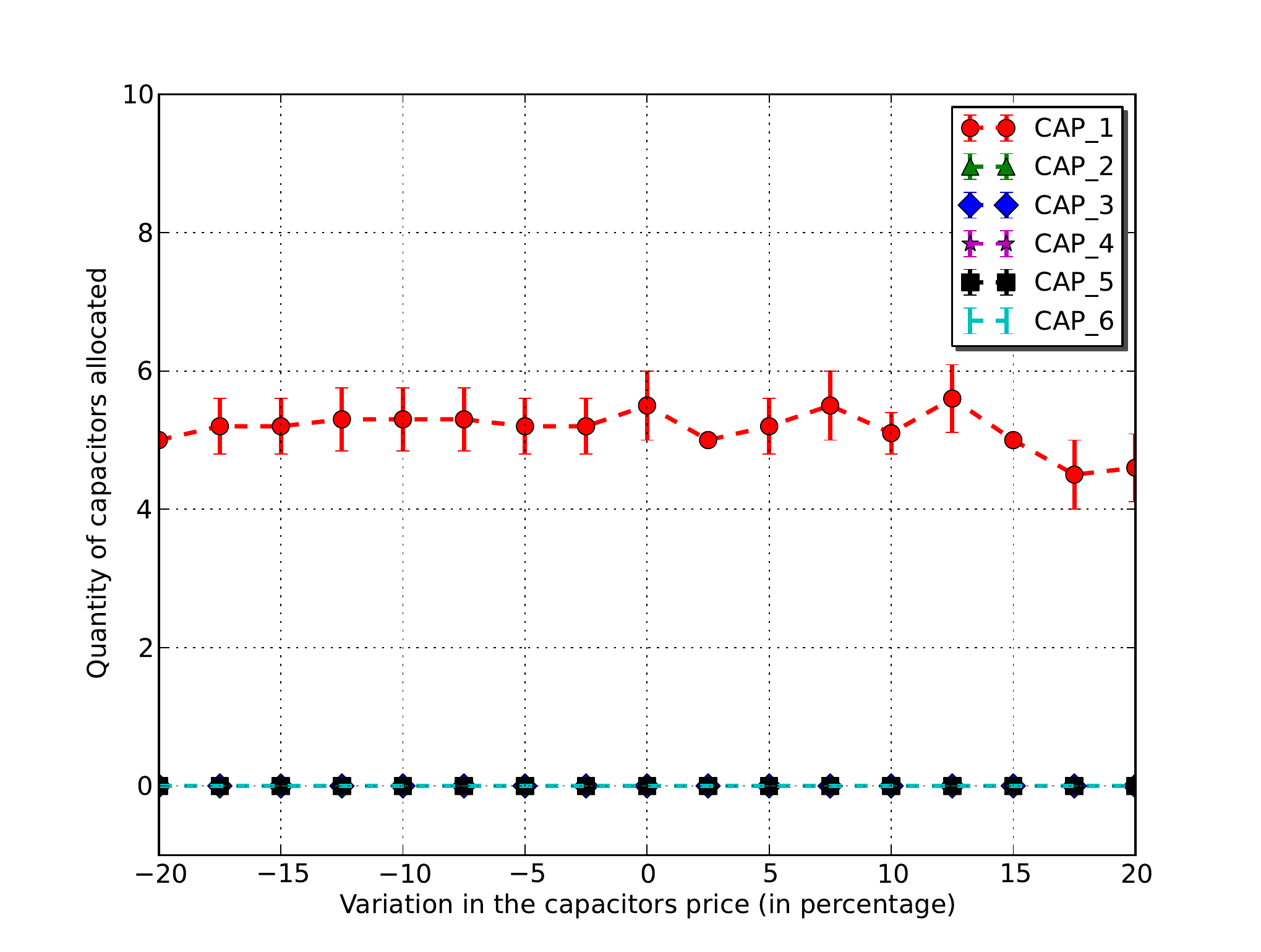}}
  \caption{Quantity of each type of capacitors allocated by the algorithms varying the capacitors cost.}
  \label{fig:type_caps_allocated_varcaps}
\end{figure*}

From now we will focus on the behavior of the algorithms when capacitor prices are changed. It is expected that the number of 
allocated capacitors will reduce as the capacitor price increases. Figures (\ref{fig:eo_qte_varcaps}), 
(\ref{fig:stg1_qte_varcaps}), and (\ref{fig:stg2_qte_varcaps}) shows the mean and standard deviation of 
the number of capacitors allocated by each type. We vary the capacitor prices from -20 to 20 percent, in other words, 
we construct scenarios of reduction and increase in the capacitor prices.

In contrast to what was expected, the number of allocated capacitors did not change significantly as the capacitor
prices increase. It shows signs of robustness of the algorithms, although more simulations, mainly considering real 
large-sized power distribution networks, are required.

The EO algorithm performs only one slight modification of the current solution per step. So, generally, 
more steps are necessary to reach high-quality solutions, as might be expected. However,
the final number of fitness evaluations is still smaller than the one required by the population-based memetic algorithm, even
adopting a structured population.

\begin{figure}[!ht]
	\centering
	\includegraphics[scale=0.4]{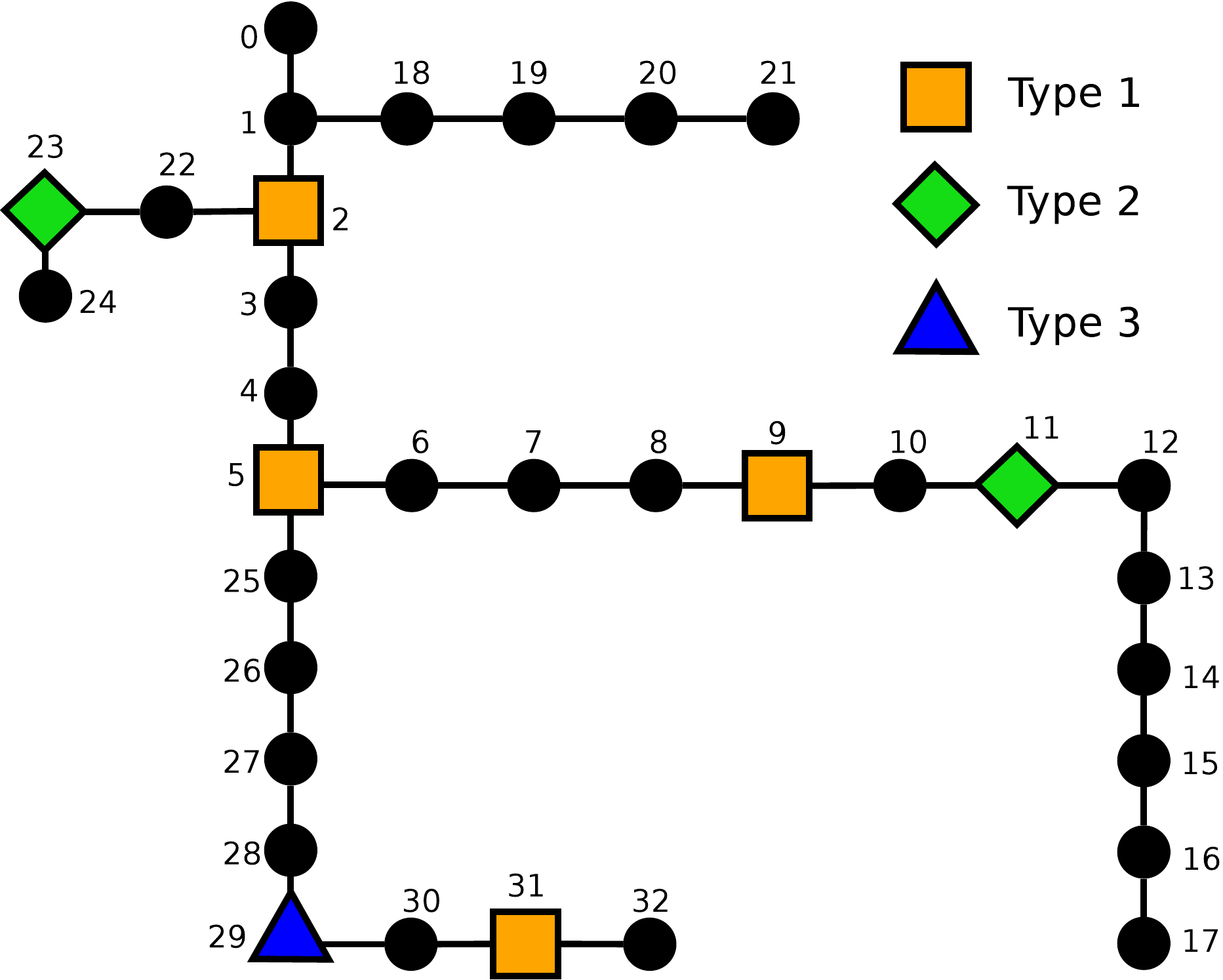}
	\caption{A solution obtained by the EO algorithm.}
	\label{fig:solucao}
\end{figure}

A solution returned by the EO algorithm is illustrated in Figure~\ref{fig:solucao}. 

\section{Conclusion}
\label{sec:conclusao}

It is known that the proper installation of capacitors brings some benefits 
to power distribution networks, such as maintenance of network stability
and reduction of loss due to reactive currents. However, this procedure 
could allocate capacitors which become harmonic sources in the network. 

In this work an Extremal Optimization based algorithm was proposed to tackle 
the optimal capacitor placement problem, besides including, in a parsimonious way,
parallel resonance constraints. 

The results showed that this approach reaches better performance when compared with 
another one which does not take resonance constraint into account, and attempts 
to repair candidate solutions by means of some post-optimization procedure. Assuming that
the optimal solution of the unconstrained problem was found, trying to make it
feasible, using local information, does not guarantee its optimality in the
constrained problem.

It should also be interesting in future works to study new ways to generate neighbors, 
maybe using some heuristics to create good solutions. Another investigation is
the behavior of our methodology when dealing with real large-sized power distribution networks,
in order to account for its scalability. As mentioned earlier, it is possible to analyze
the balance between cost of losses and investment made in capacitors from a multiobjective perspective,
this will also be the focus of future works.

\section{Acknowledgements}

We would like to thank CNPq for the financial support for this research.
The first author also like to thank Salom\~ao S. Madeiro for providing
very helpful comments.

\bibliographystyle{natbib}
\bibliography{referencias}

\end{document}